\documentclass[11pt]{amsart}
\usepackage[dvips]{graphicx}
\usepackage{amsmath,amsthm}
\usepackage{amssymb,latexsym}
\usepackage[a4paper]{geometry}
\usepackage{url}
\newtheorem{thm}{Theorem}[section]
\newtheorem{lem}[thm]{Lemma}

\newtheorem{rem}[thm]{Remark}
\newtheorem{Cor}[thm]{Corollary}

\newcommand\R{{\mathbb R}}

\newcommand\B{\overset{\circ}{B}}
\renewcommand\S{{\mathbb S}}

\newcommand\ind{{\mathbf 1}}
\newcommand\E{{\mathbb E}}
\begin{document}
\date{\today}

\author{Anne ESTRADE}
\address{Anne ESTRADE, MAP5 UMR CNRS 8145, Universit\'e Paris Descartes,
45 rue des Saints-P\`eres,  F-75270 PARIS cedex 06 \& ANR-05-BLAN-0017}
\email{anne.estrade\@@{}parisdescartes.fr}
\author{Jacques ISTAS} 
\address{Jacques ISTAS, Laboratoire Jean Kuntzmann, Universit\'e de Grenoble et CNRS, F-38041 GRENOBLE cedex 9}
\email{Jacques.Istas\@@{}imag.fr}

\title{Ball throwing on spheres}

\begin{abstract}
Ball throwing on Euclidean spaces has been considered for a
while. A suitable renormalization leads to a fractional Brownian
motion as limit object. In this paper we investigate ball throwing on spheres. A 
different behavior is exhibited: we still get a Gaussian limit but which is no
longer a fractional Brownian motion. However the limit is locally
self-similar when the self-similarity index $H$ is less than $1/2$. 
\end{abstract}

\maketitle

\section{Introduction}
Random balls models have been studied for a long time and are known as germ-grain models,
shot-noise,  or  micropulses. The common point of those
models consists in throwing balls that eventually overlap
at random in an $n$-dimensional space. Many random phenomena can  be modelized through this procedure and many
application  fields are  concerned:  Internet traffic  in  dimension one,
communication network or imaging in dimension two, biology or materials sciences in
dimension  three.  A  pioneer  work is  due  to  Wicksell
\cite{Wicksell}  with the  study  of  corpuscles.  The  
literature on germ-grain models deals with two main axes. Either
the research focuses on the geometrical or morphological aspect of the union of
random balls (see \cite{Se} or \cite{SKM} and references therein), or it is interested in the number of balls covering each
point.  This second point of view is currently known as shot-noise or
spot-noise (see \cite{Da} for existence). 
In dimension three, the  shot-noise process is a natural  candidate for modeling
porous  or  granular  media,  and  more  generally  heterogeneous  media  with
irregularities at any scale. The idea is to build
a  microscopic  model  which   yields  a  macroscopic  field  with
self-similar properties. The same idea is expected in dimension
one for Internet traffic for instance \cite{WT}. 
A usual way for catching self-similarity is  to deal
with scaling limits. Roughly speaking, the balls are dilated with a scaling parameter
$\lambda$ and one lets $\lambda$ go either to 0 or to infinity. 
We  quote  for  instance  \cite{CGM}  and  \cite{KLN}  for  the  case  $\lambda
\rightarrow 0^+$, \cite{HS} and \cite{BE} for  the  case  $\lambda
\rightarrow  +\infty$,   \cite{CT}  and   \cite{BEK}  where  both   cases  are
considered.

In the present paper, we follow a procedure which is similar to \cite{BE} and \cite{BEK}. 
Let us describe it precisely. A collection of random
balls in $\mathbb R^n$ whose centers and radii are  chosen
according to a random Poisson measure on $\mathbb R^n \times \mathbb
R^+$ is considered.  The Poisson intensity is prescribed as follows
\begin{eqnarray*}
\nu({\rm d}x,{\rm d}r) & = & r^{-n-1+2H}\,{\rm d}x\, {\rm d}r~,
\end{eqnarray*}
for some real parameter $H$. Since the Lebesgue measure ${\rm d}x$ is invariant with respect to isometry,
so is the random balls model, and so will be any (eventual) limit.
As the distribution of the radii follows a homogeneous
distribution, a self-similar scaling limit may be expected. Indeed, with additional technical
conditions, the scaling limits of such random balls models are isometry invariant
self-similar Gaussian fields. The self-similarity index depends on the parameter
$H$. When $0<H<1/2$, the Gaussian field is nothing but the
well-known fractional Brownian motion \cite{KOL,MANNESS,SATA}.

Manifold indexed fields that share properties with Euclidean self-similar fields
have been and still are extensively studied 
(e.g. \cite{FAR,GANG,Istas,ISTASMETRIC,IL,LEVY,TAK,TAKKUUR}). 
In this paper, we wonder what happens when balls are thrown onto a sphere, and no
longer onto a Euclidan space. More precisely, is there a scaling limit
of random balls models, and, when it exists, is this scaling limit a
fractional Brownian field indexed by the surface for $0<H<1/2$?

The random field is still
obtained by throwing overlapping balls in a Poissonian way. The
Poisson intensity is chosen as follows
\begin{eqnarray*}
\nu({\rm d}x,{\rm d}r) & = & f(r)\, \sigma({\rm d}x)\,{\rm d}r.
\end{eqnarray*}
The Lebesgue measure ${\rm d}x$ has been replaced by the surface
measure $\sigma({\rm d}x)$. The function $f$, that manages the
distribution of the radii, is still equivalent to $ r^{-n-1+2H}$, at
least for small $r$, where $n$ stands for the surface dimension.
It turns out that the results are completely different in the two cases
(Euclidean, spherical).
In the spherical case, there is a Gaussian scaling limit
  for any $H$. But it is no longer a fractional Brownian field, as defined by
\cite{Istas}. We then investigate the local behavior, in the tangent
bundle, of this scaling limit in the spirit of local self-similarity
\cite{BEJARO,DOB,IL}. It is locally 
asymptotically  self-similar, with a Euclidean fractional Brownian field as
tangent field. Our microscopic model has led to a local self-similar
macroscopic model.

The paper is organized as follows. In Section 2, the spherical model
is introduced   and we prove the existence of a scaling limit.  In Section 3, we study the locally self-similar
property of the asymptotic field. Section 4 is devoted to a comparative analysis between Euclidean and
spherical cases.  Eventually, some technical
computations are presented in the Appendix.


\section{Scaling limit}

We work on $\S_n$ the $n$-dimensional unit sphere, $n\ge 1$:
\begin{eqnarray*}
\S_n &= & \{ (x_i)_{1\le i \le n+1} \in \mathbb R^{n+1}\;; \; \sum_{1\le i \le n+1} x_i^2 =1 \} .
\end{eqnarray*}

\subsection{Spherical caps}
For $x,y\in \S_n$, let $d(x,y)$ denote the distance between  $x$ and
$y$ on $\S_n$, i.e. the non-oriented angle
between $Ox$ and $Oy$ where $O$ denotes the origin of $\R^{n+1}$. For 
$r \geq 0$, $B(x,r)$ denotes the closed ball  on $S$ centered
at $x$ with radius $r$:
\begin{eqnarray*} 
B(x,r) &= & \{y \in \S_n \;; \; d(x,y) \leq r \}.
\end{eqnarray*}
Let us notice that for $r< \pi$, $B(x,r)$ is a spherical cap on the unit sphere $\S_n$, 
centered at $x$ with opening angle $r$ and that
for $r \geq \pi$, $B(x,r)=\S_n$.  \\
Denoting $\sigma({\rm d}x)$  the   surface
measure on $\S_n$, we prescribe $\phi(r)$ as the surface of any ball on $\S_n$ with
radius $r$,
$$\phi(r):=\sigma(B(x,r))~,~x\in \S_n,r\ge 0~.$$
We also introduce the following function defined  
for $z$  and $z'$, two  points in $\S_n$  and $r\in \R^+$, 
\begin{equation} \label{defPsi}
\Psi(z,z',r):=\int_{\S_n} ~\ind_{d(x,z)<r}~\ind_{d(x,z')<r}~\sigma({\rm d}x)~.
\end{equation}
Actually $\Psi(z,z',r)$ denotes the surface measure of the set of all points in $\S_n$ that
belong to  both balls $B(z,r)$  and $B(z',r)$. Clearly $\Psi(z,z',r)$  only depends  on the
distance    $d(z,z')$    between    $z$    and    $z'$.    We      write
\begin{equation} \label{defpsi}
\psi(d(z,z'),r)=\Psi(z,z',r)
\end{equation} 
and note that it satisfies the following: $\forall (u,r)\in [0,\pi]\times \R^+$ 

$\bullet$  $0\le \psi(u,r) \le \sigma(\S_n)\wedge \phi(r)$

$\bullet$ if $r<u/2$ then $ \psi(u,r) = 0$ and if $r>\pi$ then $ \psi(u,r) = \sigma(\S_n)$

$\bullet$ $\psi(0,r)=\phi(r)\sim c\,r^{n}$ as $r\rightarrow 0^+~.$ 
\\~\\
In what follows we consider a family
of balls $B(X_j,R_j)$ generated at random, following a strategy described in the next section.

\subsection{Poisson point process} 
We consider a Poisson point
process $(X_j,R_j)_j$ in $\S_n\times \R^+$, or equivalently $N({\rm
d}x,{\rm d}r)$ a Poisson random measure on $\S_n\times \R^+$, with intensity
$$\nu({\rm d}x,{\rm d}r)=f(r)\,\sigma({\rm d}x)\,{\rm d}r$$
where  $f$  satisfies  the
following assumption ${\mathbf A(H)}$ for some $H>0$:

$\bullet$ $supp(f)\subset [0,\pi)$

$\bullet$ $f$ is bounded on any compact subset of $(0,\pi)$

$\bullet$ $f(r)\sim r^{-n-1+2H}~ \mbox{ as } ~r\rightarrow  0^{+}$.
\\~\\
Remarks :\\
1) The first condition ensures that no balls of radius $R_j$ on the sphere self-intersect. \\
2) Since $\phi(r)\sim c\,r^{n}~,~r\rightarrow 0^{+}$, the last condition implies that $\int_{\R^+}\phi(r)f(r)\,{\rm d}r<+\infty$, which means
that the mean surface, with respect to $f$,  of the balls $B(X_j,R_j)$ is finite.

\subsection{Random  field} 
  Let $\mathcal{M}$ denote the space of
signed measures $\mu$ on $\S_n$ with finite total variation
$|\mu|(\S_n)$, with $|\mu|$ the total variation measure of $\mu$. For
any $\mu \in \mathcal{M}$, we define 
\begin{equation}\label{Xmu}
X(\mu)=\int_{\S_n\times \R^+}\mu(B(x,r))~N({\rm d}x,{\rm d}r)~.
\end{equation}
Note that the stochastic integral in \eqref{Xmu} is well defined since 
\begin{eqnarray*}
\int_{\S_n\times \R^+}|\mu(B(x,r))|~f(r)\sigma({\rm d}x){\rm d}r 
&\le& \int_{\S_n}\int_{\S_n}  \int_{\R^+}\ind_{d(x,y)<r}\,f(r)\,\sigma({\rm d}x)\,|\mu|({\rm d}y)\,{\rm d}r\\
&&~=~ |\mu|(\S_n) \, \left( \int_{\R^+}\phi(r)f(r){\rm d}r \right)~ <+\infty~~. 
\end{eqnarray*}
In the particular case where $\mu$ is a Dirac measure $\delta_z$ for some point $z\in \S_n$ we simply denote 
\begin{equation}\label{Xz}
X(z)=X(\delta_z)=\int_{\S_n\times \R^+}\ind_{B(x,r)}(z)~N({\rm d}x,{\rm d}r)~.
\end{equation}
The pointwise field $\{X(z);z\in \S_n\}$ corresponds to the number of random balls $(X_j,R_j)$ covering each point of $\S_n$. 
Each random variable $X(z)$ has a Poisson distribution with mean $\int_{\R^+}\phi(r)f(r)\,{\rm d}r $.\\
Furthermore for any $\mu \in \mathcal{M}$,
$$\E(X(\mu))=  \mu(\S_n)\,\left(\int_{\R^+}\phi(r)f(r)\,{\rm d}r \right)$$
and 
$$Var(X(\mu))=\int_{\S_n\times  \R^+}\mu(B(x,r))^2  f(r)\,\sigma({\rm  d}x)\,{\rm
  d}r ~~\in (0,+\infty]~.$$

\subsection{Key lemma}
For $H>0$, we would like to compute  the integral 
$$\int_{\S_n\times  \R^+}\mu(B(x,r))^2 \,r^{-n-1+2H}~\sigma({\rm
  d}x)\,{\rm  d}r$$
which is  a candidate  for the  variance of  an eventually
  scaling  limit. We first introduce $\mathcal{M}^H$ the set of measures for which the above integral does
converge:
$$
\mathcal{M}^H=\mathcal{M} \mbox{ if } 2H<n~;~\mathcal{M}^H=\{ \mu \in \mathcal{M}~;~\mu(\S_n)=0 \} \mbox{ if } 2H>n~.
$$
The following lemma deals with the function $\psi$ prescribed by \eqref{defpsi}.

\begin{lem}\label{keylem} 
Let $H>0$ with $2H \ne n$. We introduce   
$$\psi^{(H)}=\psi~\mbox{if}~0<2H<n~;~\psi^{(H)}=\psi                          -
\sigma(\S_n)~\mbox{if}~2H>n~.$$
Then for all $u\in [0,\pi]$,
$$ \int_{\R^+}|\psi^{(H)}(u,r)|\,r^{-n-1+2H}{\rm d}r < +\infty\,.$$
Furthermore, denoting 
\begin{eqnarray}
\label{def_KH}
K_H(u) & = & \int_{\R^+} \psi^{(H)}(u,r)r^{-n-1+2H}{\rm d}r
\end{eqnarray} for any $u$ in $[0,\pi]$, we have 
 for all $ \mu \in {\mathcal M}^H$,
$$ 0\le \int_{\S_n\times \R^+}\mu(B(x,r))^2 r^{-n-1+2H}\,\sigma({\rm d}x)\,{\rm d}r
= \int_{\S_n\times \S_n} K_H(d(z,z'))\,\mu({\rm d}z)\mu({\rm d}z') < +\infty. $$
\end{lem}
~\\

\begin{rem}
~\\
1) For  $x,y$  in $\S_n$,  the
difference of Dirac measures $\delta_x-\delta_y$ belongs to $\mathcal{M}^H$ for any $H$.\\
2) In the case $2H>n$, since any $ \mu \in {\mathcal M}^H$ is centered, the rhs integral is not changed
when a constant is added to the kernel $K_H$.\\
3) This lemma proves that the kernel $K_H$ defines a covariance function on ${\mathcal M}^H$.
\end{rem}

\begin{proof} Using the properties of $\psi$,  we get in the case  $0<2H<n$,
\begin{eqnarray*}
0\le \int_{\R^+}\psi(u,r)\,r^{-n-1+2H}{\rm d}r &\le & 
\int_{(0,\pi)}\phi(r)\,r^{-n-1+2H}{\rm d}r + \sigma(\S_n) \int_{(\pi,\infty)}r^{-n-1+2H}{\rm d}r \\
&< &  +\infty~.
\end{eqnarray*}
In the same vein, in the case $2H>n$, we get
\begin{eqnarray*}
0\le \int_{\R^+}\left( \sigma(\S_n)-\psi(u,r)\right) \,r^{-n-1+2H}{\rm d}r &\le & 
\sigma(\S_n) \int_{(0,\pi)}\,r^{-n-1+2H}{\rm d}r  \\
&< &  +\infty~.
\end{eqnarray*}
We have just established that there exists a finite constant $C_H$ such that  
\begin{equation} \label{inegpsi}
\forall u\in [0,\pi]~,~ \int_{\R^+}|\psi^{(H)}(u,r)|\,r^{-n-1+2H}{\rm d}r\le C_H \,.
\end{equation}
The first statement is proved.
\\~\\ 
Let us denote for $ \mu \in {\mathcal M}^H$
$$I_H(\mu)=\int_{\S_n\times  \R^+}\mu(B(x,r))^2  r^{-n-1+2H}\,\sigma({\rm
  d}x)\,{\rm d}r$$
and start with proving that $I_H(\mu)$ is finite. 
We will essentially use Fubini's Theorem in the following lines. 
\begin{eqnarray*}
I_H(\mu)
&=&\int_{\R^+} \left( \int_{\S_n}\mu(B(x,r))^2\, \sigma({\rm d}x) \right) r^{-n-1+2H}{\rm
  d}r\\
&=&\int_{\R^+} \left( \int_{\S_n \times \S_n} \Psi(z,z',r)\,
\mu({\rm d}z)\mu({\rm d}z')\right) r^{-n-1+2H}{\rm d}r~.
\end{eqnarray*}
Since $  \mu\in {\mathcal M}^H$ is
centered in the case $2H>n$, one can change $\psi$ into $\psi^{(H)}$ within the previous integral. Then
\begin{eqnarray*}
I_H(\mu)
&\le& \int_{\R^+} \left( \int_{\S_n \times \S_n} |\psi^{(H)}(d(z,z'),r)|\,
|\mu|({\rm d}z)|\mu|({\rm d}z')\right) r^{-n-1+2H}{\rm d}r \\
&\le & \int_{\S_n \times   \S_n}\left(\int_{\R^+}|\psi^{(H)}(d(z,z'),r)|\,r^{-n-1+2H}{\rm
    d}r\right)|\mu|({\rm d}z)|\mu|({\rm d}z')\\
&\le & C_H ~|\mu|(\S_n)^2 <+\infty~.
\end{eqnarray*}
~\\
Following the same lines (except for the last one) without the ``$|~|$'' allows the computation of
$I_H(\mu)$ and concludes the proof.
\end{proof}
~\\
An explicit  value for  the kernel $K_H$  is available starting  from its definition. The point is  to compute $\psi^{(H)}$. 
We give in the Appendix a recurrence formula for $\psi^{(H)}$, based on the dimension $n$ of the unit
sphere $\S_n$ (see Lemma \ref{recpsi}).

\subsection{Scaling}
Let  $\rho >0$ and  $\lambda$ be  any positive  function on  $(0,+\infty)$. We
consider  the  scaled Poisson  measure  $N_\rho$  obtained  from the  original
Poisson  measure $N$  by  taking the  image  under the  map $(x,r)\in  S\times
\R^+~\mapsto~(x,\rho     r)$    and     multiplying     the    intensity     by
$\lambda(\rho)$. Hence $N_\rho$ is still a Poisson random measure with intensity
$$\nu_\rho({\rm  d}x,{\rm  d}r)=   \lambda(\rho)  \rho^{-1}  f(\rho^{-1}r)   \,\sigma({\rm  d}x)\,{\rm
  d}r~.$$
We also introduce the scaled random field $X_\rho $ defined on $\mathcal{M}$ by
\begin{equation}\label{Xmurho}
X_\rho(\mu)=\int_{\S_n\times \R^+}\mu(B(x,r))~N_\rho({\rm d}x,{\rm d}r)~.
\end{equation}

\begin{thm}\label{scaling} Let $H>0$ with $2H\ne n$ and let $f$ satisfy ${\mathbf A}(H)$.
For all positive functions $\lambda$ such that $\lambda(\rho)\rho^{n-2H}\underset
{\rho\rightarrow +\infty}{\longrightarrow}+\infty$,
the limit 
$$\left\{
\frac{X_{\rho}(\mu)-\E(X_{\rho}(\mu))}{\sqrt{\lambda(\rho)\rho^{n-2H}}};\mu\in {\mathcal M}^H\right\}\underset{\rho\rightarrow
+\infty}{\overset{fdd}{\longrightarrow}}\left\{W_H(\mu);\mu\in {\mathcal M}^H\right\}
$$ holds in the sense
of finite dimensional distributions of the random functionals. Here $W_H$ is the centered Gaussian
 random linear functional on ${\mathcal M}^H$ with
\begin{equation}\label{covarianceW}
\mbox{Cov}\left(W_H(\mu),W_H(\nu)\right)
=\int_{\S_n\times \S_n} K_H(d(z,z'))\,\mu({\rm d}z)\nu({\rm d}z'),
\end{equation}
where  $K_H$ is the kernel introduced in Lemma \ref{keylem}.
\end{thm}

The theorem can be rephrased in term of the pointwise field
$\{X(z);z\in \S_n\}$ defined in (\ref{Xz}).

\begin{Cor} \label{cor-scaling}
Let $H>0$ with $2H\ne n$ and let $f$ satisfy ${\mathbf A}(H)$.
For all positive functions $\lambda$ such that $\lambda(\rho)\rho^{n-2H}\underset
{\rho\rightarrow +\infty}{\longrightarrow}+\infty$,\\
$\bullet$ if $0<2H<n$ then  
$$\left\{
\frac{X_{\rho}(z)-\E(X_{\rho}(z))}{\sqrt{\lambda(\rho)\rho^{n-2H}}};z\in
\S_n \right\}
\underset{\rho\rightarrow +\infty}{\overset{fdd}{\longrightarrow}}\{W_H(z);z\in \S_n\}
$$
where $W_H$ is the centered Gaussian  random field on $\S_n$ with
$$\mbox{Cov}\left(W_H(z),W_H(z')\right)= K_H(d(z,z'))~.$$
$\bullet$ if $2H>n$ then for any fixed point $z_0\in \S_n$,  
$$\left\{
\frac{X_{\rho}(z)-X_{\rho}(z_0)}{\sqrt{\lambda(\rho)\rho^{n-2H}}};z\in
\S_n \right\}
\underset{\rho\rightarrow
+\infty}{\overset{fdd}{\longrightarrow}} \{W_{H,z_0}(z);z\in \S_n\}
$$ 
where $W_{H,z_0}$ is the centered Gaussian  random field on $\S_n$ with
$$\mbox{Cov}\left(W_{H,z_0}(z),W_{H,z_0}(z')\right)= K_H(d(z,z'))-K_H(d(z,z_0))-K_H(d(z',z_0))+K_H(0)~.$$
\end{Cor}

\medskip

\noindent \begin{proof} of Theorem \ref{scaling}.\\
Let us denote $n(\rho):=\sqrt{\lambda(\rho)\rho^{n-2H}}$. The 
characteristic function of the normalized field $\left({X_{\rho}(.)-\E(X_{\rho}(.))}\right)
/{n(\rho)}$ is then given by
\begin{equation}\label{charact}
\mathbb{E}\left(\exp\left({i\frac{X_{\rho}(\mu)-\E(X_{\rho}(\mu))}{n(\rho)}}\right)\right)
=\exp\left(\int_{\S_n\times  \R^+}G_\rho(x,r)\,   {\rm  d}r  \sigma({\rm  d}x)
\right)
\end{equation}
where 
\begin{equation}\label{Grho}
G_\rho(x,r)=\left( e^{i\frac{\mu(B(x,r))}{n(\rho)}}-1-i\frac{\mu(B(x,r))}{n(\rho)}
\right)\, \lambda(\rho)\rho^{-1} f(\rho^{-1}r)~.
\end{equation}
We  will  make  use of  Lebesgue's  Theorem  in  order  to  get the  limit  of
$\int_{\S_n\times  \R^+}G_\rho(x,r)\,  {\rm  d}r  \sigma({\rm d}x)$  as  $\rho\rightarrow
+\infty$. 

On the one hand, $n(\rho)$ tends to $+\infty$ 
so that $\left(e^{i\frac{\mu(B(x,r))}{n(\rho)}}-1-i\frac{\mu(B(x,r))}{n(\rho)}\right)$  behaves  like
$-\frac{1}{2}\left( \frac{\mu(B(x,r))}{n(\rho)}\right)^2$.  Together with the
assumption ${\mathbf A}(H)$, it yields the following asymptotic. For all
$(x,r)\in \S_n\times  \R^+$,
\begin{equation} \label{limGrho}
G_\rho(x,r)~\underset{\rho\rightarrow        +\infty}{\longrightarrow}       ~
-\frac{1}{2}\,\mu(B(x,r))^2 \,r^{-n-1+2H}~.
\end{equation}

On the  other hand,  since $\frac{|\mu|(B(x,r))}{n(\rho)}\le |\mu|(\S_n)$ for $\rho$ large enough,  we note
that there exists some positive constant $K$ such that for all $x,r,\rho$,
$$\left|e^{i\frac{\mu(B(x,r))}{n(\rho)}}-1-i\frac{\mu(B(x,r))}{n(\rho)}
\right|\le K\, \left( \frac{\mu(B(x,r))}{n(\rho)}\right)^2~.$$
Therefore
\begin{eqnarray*}
|G_\rho(x,r)| & \le & K \mu(B(x,r))^2 \,\rho^{-n-1+2H}\,f(\rho^{-1}r).
\end{eqnarray*}
There exists $C>0$ such that for all $r\in \R^+$~,~$f(r)\le  Cr^{-n-1+2H}$. Then 
we
get
\begin{equation} \label{majoGrho}
|G_\rho(x,r)|~\le ~ KC\,\mu(B(x,r))^2 \,r^{-n-1+2H}~
\end{equation}
where  the  right  hand  side   is  integrable  on  $\S_n\times  \R^+$  by  Lemma
\ref{keylem}. 

Applying Lebesgue's Theorem yields
\begin{eqnarray*}
\int_{\S_n\times \R^+}G_\rho(x,r)\, \sigma({\rm d}x)\,{\rm d}r
~\underset{\rho\rightarrow +\infty}{\longrightarrow} 
&&  -\frac{1}{2}\, \int_{\S_n\times  \R^+}\mu(B(x,r))^2  \,r^{-n-1+2H}\,\sigma({\rm d}x)\,{\rm d}r\\
&=& -\frac{1}{2}\,\int_{
\S_n\times \S_n} K_H(d(z,z'))\,\mu({\rm d}z)\mu({\rm d}z')~.
\end{eqnarray*}

Hence $(X_{\rho}(\mu)-\E(X_{\rho}(\mu)))/n(\rho)$ converges in
distribution to the centered Gaussian random variable $W(\mu)$ whose
variance is equal to
$$\mathbb{E}\left(W(\mu)^2\right)=C\,\int_{\S_n\times \S_n} K_H(d(z,z'))\,\mu({\rm d}z)\mu({\rm d}z')~.$$
By linearity, the covariance of $W$ satisfies (\ref{covarianceW}). 

\end{proof}

\begin{rem}
~\\
1) The pointwise limit field $\{W_H(z);z\in \S_n\}$
in Corollary \ref{cor-scaling} is stationary, i.e. its distribution is  invariant under the isometry group of $\S_n$, whereas the increments of
$\{W_{H,z_0}(z);z\in \S_n\}$ are  distribution invariant under  the group of all  isometries of
$\S_n$ which keep the point $z_0$ invariant.  \\
2) When $0<H<1/2$ the Gaussian field $W_H$ does not coincide with the field introduced in
\cite{Istas} as the spherical fractional Brownian motion on $\S_n$.
\end{rem}
 Indeed, let us have a look at the case $n=1$, it is easy to obtain the following piecewise expression for
$\psi=\psi_1$: $\forall (u,r)\in [0,\pi]\times \R^+$,
\begin{eqnarray*}
\psi_1(u,r)&=&0 ~~\mbox{ for } 0\le r<u/2 \\
&=&2r-u ~~\mbox{ for } u/2\le r \le \pi-u/2 \\
&=&4r-2\pi~~ \mbox{ for } \pi-u/2 \le r \le \pi\\
&=&2\pi ~~\mbox{ for }r> \pi
\end{eqnarray*}
and to compute
\begin{eqnarray*}
K_H(u) & = &\frac{1}{H(1-2H)2^{2H}}\left( 2(2H)^{2H}-u^{2H}-(2\pi-u)^{2H} \right)~,
\end{eqnarray*}
Actually, we compute the variance of the increments of $W_H$
\begin{eqnarray*}
\mathbb E(W_H(z)-W_H(z'))^2 & = & 2 K_H(0)-2K_H(d(z,z')) \\
 & = &  \frac{2}{H(1-2H) 2^{2H}} [ d^{2H}(z,z') +(2 \pi -d(z,z'))^{2H} -(2 \pi)^{2H} ].
\end{eqnarray*}
The spherical fractional Brownain motion $B_H$, introduced in \cite{Istas}, satisfies
\begin{eqnarray*}
\mathbb E(B_H(z)-B_H(z'))^2 & = &  d^{2H}(z,z').
\end{eqnarray*}
Even up to a constant, processes $W_H$ and $B_H$ are clearly
different. The Euclidean situation is therefore different. Indeed,
\cite{BEK}, the variance of the increments of the corresponding  field
$W_H$ is proportional to $|z-z'|^{2H}$.


\section{Local self-similar behavior}
We  wonder whether  the limit  field $W_H$  obtained in  the  previous section
satisfies a local asymptotic self-similar ({\it lass}) property.
More precisely we will let a ``dilation'' of order $\varepsilon$ act on $W_H$ near a fixed point $A$ in $\S_n$ and as in \cite{IL}, up to a renormalization factor,
we look for an asymptotic behavior as $\varepsilon$ goes to $0$. An $H$-self-similar tangent field $T_H$ is expected.  
Recall that $W_H$ is defined on a subspace ${\mathcal M}^H$ of measures on $\S_n$, so that $T_H$ will be defined on a subspace of measures on the 
tangent space $\mathcal T_A\S_n$ of $\S_n$.

\subsection{Dilation}

Let us fix a point  $A$ in $\S_n$ and consider $\mathcal T_A\S_n$ the tangent space of $\S_n$ at $A$. It can be identified with 
 $\R^n$ and $A$ with the null vector of $\R^n$.\\

Let $1<\delta < \pi$. The exponential map at point $A$, denoted by $\exp$, is a diffeomorphism between the Euclidean
ball  $\{y\in  \R^n,\,\|y\|<\delta\}$  and  $\B(A,\delta)\subset  \S_n$,  where
$\|.\|$ denotes the Euclidean norm  in $\R^n$ and $\B(A,\delta)$ the open ball
with center $A$ and radius $\delta$ in $\S_n$.\\
Furthermore for all $y,y'\in \R^n$ such that $\|y\|,\|y'\|<\delta$, 
$$d(A,\exp y)=\|y\| ~\mbox{ and }~d(\exp y,\exp y')\le \|y-y'\|~.$$
We refer to \cite{He} for precisions on the exponential map.
\\~\\
Let $\tau$ be a signed measure on  $\mathbb R^n$. We define the dilated measure  $\tau_{\varepsilon}$ by 
\begin{eqnarray*}
\forall B \in {\mathcal B}(\mathbb R^n) \quad \tau_{\varepsilon}(B) & = & \tau(B/\varepsilon)
\end{eqnarray*}
and then map it by the application $\exp$, defining the measure 
$\mu_{\varepsilon}=\exp^*\tau_{\varepsilon}$ on $\B(A,\delta)$ by
\begin{eqnarray} \label{dilation}
\forall C \in {\mathcal B}(\B(A,\delta)) \quad
\mu_{\varepsilon}(C)= \exp^*\tau_{\varepsilon}(C)& =& \tau_{\varepsilon}(\exp^{-1}(C)) .
\end{eqnarray}
We then consider the measure $\mu_{\varepsilon}$ as a measure on the whole sphere 
$\S_n$ with support included in $\B(A,\delta)$.
\\~\\
At last, we  define the  dilation of  $W_H$  within a  ``neighborhood of  $A$'' by  the
following  procedure. For  any finite  measure $\tau$  on $\mathbb R^n$,  we consider
$\mu_{\varepsilon}=\exp^*\tau_{\varepsilon}$  as  defined by  \eqref{dilation}
and  compute $W_H(\mu_{\varepsilon})$.  We will  establish the  convergence in
distribution of $\varepsilon^{-H}W_H(\exp^*\tau_{\varepsilon})$ for any $\tau$
in an appropriate  space of measures on $\mathbb R^n$.  Since $W_H(\mu_{\varepsilon})$ is
Gaussian, we will focus on its variance.

\subsection{Asymptotic of the kernel $K_H$}
For $0<H<1/2$, we already mentioned that the kernel $K_H(0)-K_H(u)$ is not proportional to
$u^{2H}$.   As    a   consequence,   one    cannot   expect   $W_H$    to   be
self-similar.  Nevertheless, as we are looking for an asymptotic local
self-similarity, only the behavior of $K_H$ near zero is
relevant. Actually we will establish that, up to a constant, $K_H(0)-K_H(u)$ behaves like $u^{2H}$ when $u\rightarrow 0^+$.

\begin{lem} \label{DL de KH}
Let $0<H<1/2$. The kernel $K_H$ defined by \eqref{def_KH} satisfies 
$$ K_H(u) = K_1-K_2 u^{2H} +o(u^{2H}), ~u\rightarrow 0^+$$
where $K_1=K_H(0)$ and $K_2$ are nonnegative constants. 
\end{lem}
~\\
\begin{proof}
Let us state that the assumption $H<1/2$ implies $H<n/2$ so that in that case $K_H$ is prescribed by 
$$K_H(u)=\int_{\R^+} \psi(u,r)\,r^{-n-1+2H}{\rm d}r~,~u\in [0,\pi]~.$$
We note that $K_H(0)<+\infty$ since $\psi(0,r)\sim cr^n$ as $r\rightarrow 0^+$ and $\psi(0,r)=\sigma(\S_n)$ for $r>\pi$.
Then, subtracting $K_H(0)$ and remarking that $\psi(0,r)=\psi(u,r)=\sigma(\S_n)$
for $r>\pi$, we write
\begin{eqnarray*}
K_H(0)-K_H(u)
&=&\int_0^{\pi} (\psi(0,r)-\psi(u,r))\,r^{-n-1+2H}{\rm d}r \\
&=& \int_{0}^{\delta}(\psi(0,r)-\psi(u,r))\,r^{-n-1+2H}{\rm d}r\\ 
&&~+\int_{\delta}^{\pi}(\psi(0,r)-\psi(u,r))\,r^{-n-1+2H}{\rm d}r
\end{eqnarray*}
where we recall that $\delta \in (1,\pi)$ is such that the exponential map is a
diffeormorphism     between     $\{\|y\|<\delta\}     \subset    \R^n$     and
$\B(A,\delta)\subset \S_n$.\\
The  second  term  is  of  order  $u$,  and  therefore  is
 negligible with respect to $u^{2H}$,  since $\psi$ is clearly Lipschitz on the
compact interval $[\delta,\pi]$.\\ 
We now focus on the first term. Performing the change of variable $r\mapsto
r/u$, we write it as 
$$\int_{0}^{\delta}(\psi(0,r)-\psi(u,r))\,r^{-n-1+2H}{\rm d}r=u^{2H}\; \int_{\R^+} \Delta(u,r)\,r^{-n-1+2H}{\rm d}r~,$$
where 
$$\Delta(u,r):=\ind_{ur<\delta}~u^{-n}\,\left( \psi(u,ur)-\psi(0,ur)\right)~.$$
Their only remains to prove that $\int_{\R^+} \Delta(u,r)\,r^{-n-1+2H}{\rm d}r$
admits a finite limit $K_2$ as $u\rightarrow 0^+$.
We will use Lebesgue's Theorem and start with establishing the simple convergence of $\Delta(u,r)$ for any given $r\in \R^+$.
\\~\\
We fix a unit vector $\mathbf{v}$ in $\R^n$ and a point $A'=\exp \mathbf{v}$ in $\S_n$. We then 
consider for any $u\in (0,\delta)$, the point $A'_u:=\exp(u\mathbf{v}) \in \S_n$ located on the geodesic between $A$ and $A'$ such that $d(A,A'_u)=\|u\mathbf{v}\|=u$. 
We can then use \eqref{defPsi} and \eqref{defpsi} to write
$$\psi(u,.)=\Psi(A,A'_u,.)= \int_{\S_n}~\ind_{d(A,z)<.}~\ind_{d(A'_u,z)<.}~{\rm d}\sigma(z)$$
and 
$$\psi(0,.)=\Psi(A,A,.)=\int_{\S_n}~\ind_{d(A,z)<.}~{\rm d}\sigma(z)$$
in order to express $\Delta(u,r)$ as
\begin{eqnarray*}
\Delta(u,r)
&=& \ind_{ur<\delta}~u^{-n}\,\int_{\S_n} \ind_{d(A,z)<ur}\,\ind_{d(A'_u,z)>ur}\,{\rm d}\sigma(z)~.
\end{eqnarray*}
Since $ur<\delta$ the above integral runs on $\B(A,ur)\subset \B(A,\delta)$ and we can perform the exponential change of variable to get 
\begin{eqnarray*}
\Delta(u,r)
&=&\ind_{ur<\delta}~u^{-n}\, \int_{\R^n} \ind_{\|y\|<ur}
\,\ind_{d(\exp(u\mathbf{v}),\exp(y))>ur}\,{\rm d}\sigma(\exp(y))\\
&=& \ind_{ur<\delta}~\int_{\R^n} \ind_{\|y\|<r}\,\ind_{d(\exp(u\mathbf{v}),\exp(uy))>ur}\,\tilde{\sigma}(uy){\rm d}y~.
\end{eqnarray*}
In  the last  integral,  the image  by  $\exp$ of  the  surface measure  ${\rm
  d}\sigma(\exp(y))$ is written as $\tilde{\sigma}(y){\rm d}y$ where ${\rm d}y$ denotes the Lebesgue measure on $\R^n$. \\
We use  the fact  that $d(\exp(ux),\exp(ux'))\sim u\|x-x'\|$  as $u\rightarrow
0^+$ to get the following limit for the integrand 
$$\ind_{d(\exp(u\mathbf{v}),\exp(uy))<ur}\,\tilde{\sigma}(uy)\longrightarrow
\ind_{\|\mathbf{v}-y\|>r}\,\tilde{\sigma}(0)~.$$
Since the integrand is clearly dominated 
by 
$$\|\sigma\|_{\infty}:=\sup\{\tilde{\sigma}(y)\,,\,\|y\|\le \delta\}~,$$
Lebesgue's Theorem yields for all $r\in \R^+$,
\begin{eqnarray*}
\Delta(u,r) &\longrightarrow& \tilde{\sigma}(0)\,\int_{\R^n} \ind_{\|y\|<r}\,\ind_{\|\mathbf{v}-y\|>r}\,{\rm d}y~.
\end{eqnarray*}
~\\
We recall  that $d(\exp  x,\exp x')\le \|x-x'\|$  for all $x,x'\in  \R^n$ with
norm less than $\delta$. Therefore for all $u$,
$$\Delta(u,r) \leq \|\sigma\|_{\infty}~\int_{\R^n} \ind_{\|y\|<r}\,\ind_{\|\mathbf{v}-y\|>r}\,{\rm d}y$$ 
where the right hand side belongs to $L^1(\R^+,r^{-n-1+2H}\,{\rm d}r)$
(see \cite{BE} Lemma A.2).
\\~\\
Using Lebesgue's Theorem for the last time, we obtain
$$\int_{\R^+} \Delta(u,r)\,r^{-n-1+2H}{\rm d}r~\underset{u\rightarrow
0^+}{\longrightarrow}~K_2$$
where
$$K_2=\tilde{\sigma}(0)\,\int_{\R^n \times \R^+} \ind_{\|y\|<r}\,\ind_{\|\mathbf{v}-y\|>r}\,r^{-n-1+2H}\,{\rm d}y\,{\rm d}r~\in (0,+\infty)~.$$
\end{proof}

Let us remark that the proof makes it clear that the case $H>1/2$ is dramatically
different.  The  kernel $K_H(0)-K-H(u)$ behaves like  $u$ near zero and  looses its $2H$
power. 
\\~\\

\subsection{Main result }
~\\
Let  $0<H<1/2$. We  consider  the  following space  of  measures on  $\mathcal
T_A\S_n \cong \R^n$
\begin{eqnarray*}
\mathfrak{M}^H &=& \{\mbox{measures $\tau$  on $\mathbb R^n$ with finite total
  variation such that} \\
&&\tau(\mathbb R^n)=0 \mbox{ and }
\int_{\mathbb R^n  \times \mathbb R^n}\|x-x'\|^{2H}|\tau|({\rm d}x)|\tau|({\rm d}x')<+\infty \}~.
\end{eqnarray*}
For  any  measure  $\tau  \in  \mathfrak{M}^H$, we  compute  the  variance  of
$W_H(\mu_{\varepsilon})$
where 
$\mu_{\varepsilon}=\exp^*\tau_{\varepsilon}$ is defined by \eqref{dilation}.
\\~\\
By   Lemma   \ref{keylem},   since
$\mu_{\varepsilon}$ belongs to ${\mathcal M}={\mathcal M}^H$ in the case $H<1/2$,
\begin{eqnarray*}
\mbox{var}(W_H(\mu_{\varepsilon})) & =& \int_{B(A,\delta)\times
B(A,\delta)} K_H(d(z,z')) 
\mu_{\varepsilon}({\rm d}z)\mu_{\varepsilon}({\rm d}z') .
\end{eqnarray*}
Performing an exponential change of variable followed by a dilation in $\R^n$,
we get
\begin{eqnarray*}
\mbox{var}(W_H(\mu_{\varepsilon})) 
& =& \int_{\R^n\times \R^n} \ind_{\|y\|<\delta}\,\ind_{\|y'\|<\delta}\,K_H(d(\exp(y),\exp(y')))
\tau_{\varepsilon}({\rm d}y)\tau_{\varepsilon}({\rm d}y') \\
 & =& \int_{\R^n \times \R^n}\ind_{\|x\|<\delta/\varepsilon}\,\ind_{\|x'\|<\delta/\varepsilon}\,
K_H(d(\exp(\varepsilon x),\exp(\varepsilon x'))) 
\tau({\rm d}x)\tau({\rm d}x').
\end{eqnarray*}
Denoting  $\widetilde{K_H}(u)=K_H(u)-K_H(0)$,
\begin{eqnarray*}
\mbox{var}(W_H(\mu_{\varepsilon})) 
& =& \int_{\R^n \times \R^n}\ind_{\|x\|<\delta/\varepsilon}\,\ind_{\|x'\|<\delta/\varepsilon}\,
\widetilde{K_H}(d(\exp(\varepsilon x),\exp(\varepsilon x'))) 
\tau({\rm d}x)\tau({\rm d}x')\\
&&~+~ K_H(0)~\tau(\{\|x\|<\delta/\varepsilon\})^2~.
\end{eqnarray*}
Let us admit for a while that
\begin{eqnarray} \label{limtau}
\frac{\tau(\{\|x\|<\delta/\varepsilon\})^2}{\varepsilon^{2H}} & \underset{\varepsilon\rightarrow
0^+}{\longrightarrow} & 0~.
\end{eqnarray}
Then,   applying  Lebesgue's   Theorem  with   the  convergence   argument  on
$\widetilde{K_H}$ obtained in Lemma \ref{DL de KH}, 
yields 
\begin{eqnarray} \label{limWH}
\frac{\mbox{var}(W_H(\mu_{\varepsilon}))}{\varepsilon^{2H}} & \underset{\varepsilon\rightarrow
0^+}{\longrightarrow} &
-K_2 \int_{\mathbb R^n \times \mathbb R^n} \|x-x'\|^{2H}\tau({\rm d}x)\tau({\rm d}x').
\end{eqnarray}
~\\
Let us now establish \eqref{limtau} where we recall that $\tau$ is any measure in
$\mathfrak{M}^H$. In particular, the total mass of $\tau$ is zero so that 
\begin{eqnarray*} 
\frac{\tau(\{\|x\|<\delta/\varepsilon\})}{\varepsilon^{H}} 
&=& -\,\frac{\tau(\{\|x\|>\delta/\varepsilon\})}{\varepsilon^{H}} \\
&=& -\,\int_{\R^n}\varepsilon^{-H}\,\ind_{\|x\|>\delta/\varepsilon}\,\tau({\rm d}x)~.
\end{eqnarray*}
For        any         fixed        $x        \in         \R^n$,        $
\varepsilon^{-H}\,\ind_{\|x\|>\delta/\varepsilon}$ is zero when $\varepsilon$ is
small enough. Moreover  $\varepsilon^{-H}\,\ind_{\|x\|>\delta/\varepsilon}$     is     dominated     by
$\delta^{-H}\,\|x\|^H$  which belongs  to  $L^1(\R^n,|\tau|({\rm d}x))$  since
$\tau$ belongs to $\mathfrak{M}^H$.
Lebesgue's Theorem applies once more.
\\~\\
We deduce from asymptotic \eqref{limWH} the following theorem.

\begin{thm}
\label{sphere_lass} Let $0<H<1/2$. The limit 
$$
\frac{W_H(\exp^*\tau_{\varepsilon})}{\varepsilon^{H}}\underset{\varepsilon\rightarrow
0^+}{\overset{fdd}{\longrightarrow}}T_H(\tau)
$$ holds for all $\tau\in \mathfrak{M}^H$,  in the sense
of finite dimensional distributions of the random functionals. Here $T_H$ is the centered Gaussian  random linear functional on $\mathfrak{M}^H$ with
\begin{equation}\label{covarianceT}
\mbox{Cov}\left(T_H(\tau),T_H(\tau')\right)
=-K_2\int_{\mathbb R^n \times \mathbb R^n}
 \|x-x'\|^{2H}\,\tau({\rm d}x)\tau'({\rm d}x'),
\end{equation}
\end{thm}
~\\
As for Theorem \ref{scaling}, Theorem \ref{sphere_lass} can be
rephrased in terms of  pointwise fields. Indeed, $\delta_x - \delta_O$
belongs to $\mathfrak{M}^H$ for all $x$ in $\mathbb R^n$. Let us apply
Theorem \ref{sphere_lass} with $\tau= \delta_x - \delta_O$. Then
$T_H(\delta_x - \delta_O)$  has the covariance
\begin{eqnarray*}
\mbox{Cov}\left(T_H(\delta_x - \delta_O),T_H(\delta_{x'} - \delta_O)\right)
=K_2 (\|x\|^{2H}+ \|x'\|^{2H}-  \|x-x'\|^{2H})~,
\end{eqnarray*}
and the field $\{T_H(\delta_x - \delta_O)\,;\,x\in \R^n\}$ is a Euclidean
fractional Brownian field.

\section{Comparative analysis}

In this section, we aim to discuss the differences and the analogies between the
Euclidean and the spherical case. 
\\~\\
Let us first be concerned with the existence of a scaling limit random
field.  The variance of this limit field should be
\begin{eqnarray*}
\mathbf{V} & = & \int_{\mathbb M_n} \int_{\mathbb R^+} \mu(B(x,r))^2
\sigma({\rm d}x) r^{-n-1+2H} {\rm d}r,
\end{eqnarray*}
where $\mathbb M_n$ is the $n$-dimensional corresponding surface with its surface
measure $\sigma$. When speaking of the Euclidean case $\mathbb
M_n=\mathbb R^n$ we refer to \cite{BEK}. In the present paper, we studied the case $\mathbb M_n=\mathbb
S_n$. Moreover, in this discussion, the hyperbolic case
$\mathbb M_n = \mathbb H_n=  
\{ (x_i)_{1\le i \le n+1} \in \mathbb R^{n+1}\;; \; x_{n+1}^2-\sum_{1\le i \le n} x_i^2 =1\,,\,x_{n+1}\ge 1 \}$
is evoked.

In the Euclidean case, the random  fields are defined on the space of measures
with vanishing total mass. So let us first consider measures $\mu$ such that
$\mu(\mathbb M_n)=0$. Hence, whatever  the surface $\mathbb  M_n$, the  integral $\mathbf{V}$  involves the
integral of the surface of the 
symmetric difference between two balls of same radius $r$. As $r$ goes to infinity, three different behaviors emerge.
\begin{itemize}
\item $\mathbb M_n=\mathbb S_n$: this surface vanishes
\item $\mathbb M_n=\mathbb R^n$: the order of magnitude of this surface is $r^{n-1}$
\item $\mathbb M_n=\mathbb H_n$: the surface grows exponentially
\end{itemize}
The consequences are the following. 
\begin{itemize}
\item $\mathbb M_n=\mathbb S_n$: any positive $H$ is admissible
\item $\mathbb M_n=\mathbb R^n$: the range of admissible $H$ is $(0,1/2)$
\item $\mathbb M_n=\mathbb H_n$: no $H$ is admissible
\end{itemize}
In the Euclidean case, the restriction $\mu(\mathbb R^n)=0$ is mandatory whereas it is unnecessary in the spherical case for $H<n/2$. Indeed the integral $\mathbf{V}$ is clearly convergent. 
\\~\\
Let us now discuss the (local)  self-similarity of the limit field. Of course,
we no longer consider the hyberbolic case.
\begin{itemize}
\item $\mathbb M_n=\mathbb R^n$: dilating a ball is a homogeneous operation. Therefore, the limit
field is self-similar.
\item $\mathbb M_n=\mathbb S_n$: dilation is no longer homogeneous. Only local self-similarity can be
expected.  The natural framework of  this local self-similarity is the
tangent bundle, where the situation is Euclidean. Hence we have to come back to the restricting condition $H<1/2$.
\end{itemize}


\section*{Appendix}

\noindent {\bf Recurrence formula for the $\psi_n$'s.}

\medskip

Recall that the functions $\psi_{n}$'s are defined by \eqref{defPsi} and \eqref{defpsi}
$$
\psi_n(u,r)=\Psi_n(M,M',r)  =  \int_{\S_n}  ~\ind_{d(M,N)<r}
~\ind_{d(M',N)<r}\, {\rm d}\sigma_n(N) \; , \;(u,r)\in [0,\pi]\times \R^+\;,
$$
for any pair $(M,M')$ in $\S_n$ such that $d(M,M')=u$. Here $\sigma_n$
stands for the surface measure on $\S_n$. 

\begin{lem}
\label{recpsi}
The  family  of  functions  $\psi_n,   n  \geq  2$  satisfies  the  following
recursion:\\
 $\forall (u,r)\in [0,\pi]\times \R^+~,$
\begin{eqnarray*}
\psi_n(u,r) & = & \int_{-\sin r}^{\sin r}
(1-a^2)^{n/2} \, \psi_{n-1}\left(u,\arccos \left(\frac{\cos r}{\sqrt{1-a^2}}
\right)\right) {\rm d}a \; .
\end{eqnarray*}
\end{lem}
~\\
\begin{proof}
An arbitrary point of $\S_n$ is parameterized either in Cartesian coordinates,
$(x_i)_{1\le i\le n+1}$, or in spherical ones
$$(\phi_i)_{1\le i\le n} \in [0,\pi)^{n-1}\times [0,2\pi)$$
with
\begin{eqnarray*}
x_1 & = &  \cos \phi_1 \\
x_2 & = & \sin \phi_1 \cos \phi_2 \\
x_3 & = & \sin \phi_1 \sin \phi_2 \cos \phi_3 \\
&\cdots & \\
x_n & = & \sin \phi_1 \sin \phi_2 \ldots \sin \phi_{n-1} \cos \phi_n \\
x_{n+1} & = & \sin \phi_1 \sin \phi_2 \ldots \sin \phi_{n-1} \sin \phi_n 
\end{eqnarray*}

Let $M$  be the point
$(\phi_i)_{1\le  i\le n}=(0,\ldots,0)$.   One can  write the ball
$B_n(M,r)$ of  radius $r$,  which is  a spherical cap  on $\S_n$  with opening
angle $r$ as  follows,
\begin{eqnarray*}
B_n(M,r) & = & \{ (\phi_i)_{1\le i\le n} \in \S_n~;~\phi_1 \leq r\} 
\end{eqnarray*}
or in Cartesian coordinates
\begin{eqnarray*}
B_n(M,r) & = & \{ (x_i)_{1\le i\le n+1} \in \S_n~;~ x_1 \geq \cos r \} \; .
\end{eqnarray*}
Let $a\in (-1,1)$ and let $P_{a}$ be the hyperplane of $\R^{n+1}$ defined by $x_{n+1}=a$.
Let us consider the intersection $P_{a} \cap B_n(M,r)$. 
\begin{itemize}
\item If $1-a^2 < \cos^2 r$ then $P_{a} \cap B_n(M,r)= \emptyset$.
\item If $1-a^2 \geq \cos^2 r$ then
\begin{eqnarray*}
P_{a} \cap B_n(M,r)
&=&\{ (x_i)_{1\le i\le  n+1} \in \S_n~;~  x_1 \geq
\cos r \mbox{ and } x_{n+1}=a \}\\
&=&\{ (x_i)_{1\le i\le  n} \in \R^n~;~  x_1 \geq
\cos r \mbox{ and }\sum_{1\le i\le n}x_i^2=1-a^2\} \times \{a\}.
\end{eqnarray*}
In  other words, denoting  $\S_{n-1}(R)$  the $(n-1)$-dimensional  sphere of
radius $R$,
$$P_{a} \cap B_n(M,r)= B_{n-1,\sqrt{1-a^2}}(M(a),r(a)) \times \{a\}$$  
where $B_{n-1,\sqrt{1-a^2}}(M(a),r(a))$ is the spherical cap on $\S_{n-1}(\sqrt{1-a^2})$, centered 
at $M(a)=(\sqrt{1-a^2},0,\ldots,0)$ and with opening angle $r(a)=\arccos\left( \frac{\cos r}{\sqrt{1-a^2}}\right)$.
\end{itemize}
~\\
Let now  $M'$ be defined in spherical coordinates by $(\phi_i)_{1\le i\le n}=(u,0,\ldots,0)$, so that $d(M,M')=u$.
The intersection  $P_{a} \cap B_n(M',r)$ is the map of $P_{a} \cap B_n(M,r)$ by the rotation of angle $u$
and center $C$ in the plane $x_3=\ldots=x_{n+1}=0$. So
\begin{itemize}
\item if $1-a^2 < \cos^2 r$ then  $P_{a} \cap B_n(M',r)= \emptyset$.
\item if $1-x_0^2 \geq \cos^2 r$ then
$$P_{a} \cap B_n(M',r)= B_{n-1,\sqrt{1-a^2}}(M'(a),r(a)) \times \{a\}$$  
where the $(n-1)$-dimensional spherical cap $B_{n-1,\sqrt{1-a^2}}(M'(a),r(a))$ is now  centered
at $M'(a)=(\sqrt{1-a^2}\cos u,\sqrt{1-a^2} \sin u,0,\ldots,0)$.
\end{itemize}
~\\
We define $\psi_{n-1,R}(u,r)$ as the intersection surface  of two
spherical caps on $\S_{n-1}(R)$, whose centers are at a distance $Ru$ and with the same opening angle
$r$.\\
By homogeneity, this leads to
$$\psi_{n-1,R}(u,r)=R^{n-1}\,\psi_{n-1,1}(u,r)=R^n\,\psi_{n-1}(u,r)\,.  $$
The surface measure $\sigma_{n}$ of $\S_n$ can be written as
$${\rm d}\sigma_{n}(x_1,\ldots,x_n,a)=
\sqrt{1-a^2}\;{\rm d}\sigma_{n-1,\sqrt{1-a^2}}(x_1,\ldots,x_n)\times {\rm d}a$$
where $\sigma_{n-1,R}$ is the surface measure of $\S_{n-1}(R)$.\\
We then obtain
\begin{eqnarray*}
\psi_n(u,r) & = & \int_{-1}^1  ~\ind_{1-a^2 \geq \cos^2 r}
~\psi_{n-1,\sqrt{1-a^2}}\,(u,\arccos\left( \frac{\cos r}{\sqrt{1-a^2}}\right)) \sqrt{1-a^2} {\rm d}a \\
 & = & \int_{-\sin r}^{\sin r}(\sqrt{1-a^2})^n ~\psi_{n-1}(u,\arccos\left( \frac{\cos r}{\sqrt{1-a^2}}\right)) {\rm d}a \,, 
\end{eqnarray*} 
and Lemma \ref{recpsi} is proved.
\end{proof}


\section*{Acknowledgements}
The authors are  very greatful to the referee for  his/her careful reading and
useful comments.


\bibliographystyle{plain}

\end{document}